\documentclass[12pt]{article}

 \usepackage[utf8]{inputenc}
\usepackage[english]{babel}
\usepackage{amsfonts,amsmath,amsxtra,amsthm,latexsym}
\usepackage{amssymb} 
 \theoremstyle{definition}
 
 \theoremstyle{remark}
 \newtheorem{rem}{Remark}[section]

  \newtheorem*{trem*}{Draft remark}

 \numberwithin{equation}{section}

\DeclareMathOperator{\dom}{dom}

\def\RR{\mathbb R} \def\CC{\mathbb C}

 \def\DD{\mathbb D}

\def\MM{\mathbb{M}}

\newcommand{\x}{\mathbf{x}}

\def\al{\alpha}
\def\th{\theta}

\def\de{\delta}
\def\vep{\varepsilon}
\def\ep{\epsilon}

\def\Ga{\Gamma}
\def\om{\omega}
\def\Om{\Omega}
\def\Si{\Sigma}
\def\si{\sigma}

\def\pa{\partial}

\def\ii{\mathrm{i}}

\def\dd{\mathrm{d}}

\def\modn{\hspace{-8pt}\mod}
\def\<{\langle}
\def\>{\rangle}

\def\M{\mathcal{M}}

\def\out{\mathrm{out}}
\def\abs{\mathrm{abs}}

\def\epr{\vep_{\mathrm{r}}}
\def\epout{\vep_{\mathrm{out}}}
\def\E{\mathbf{E}}
\def\H{\mathbf{H}}
\def\n{\mathbf{n}}
\def\Z{\mathcal{Z}}

\def\C{\mathcal{C}}
\def\K{\mathcal{K}}
\def\c{\mathrm{opt}}

\begin{document}

\title{Euler-Lagrange equations for full topology optimization of the Q-factor in leaky cavities} 
\author{}
\date{}
\maketitle
  
{ \center{ \large
Matthias Eller$^\text{ a}$ and Illya M. Karabash$^\text{ b,c}$\\[3ex]
  }
}  
  
  {\small \noindent
$^{\text{a}}$ Department of Mathematics and Statistics, Georgetown University, Washington, DC 20057 \\[1mm]  
$^{\text{b}}$ Mathematical Institute, Rheinische Friedrich-Wilhelms Universität Bonn,  Endenicher  Allee 60, D-53115 Bonn, Germany \\[1mm]
 $^{\text{c}}$ Institute of Applied Mathematics and Mechanics of NAS of Ukraine,
Dobrovolskogo st. 1, Slovyans'k 84100, Ukraine\\[1mm]
E-mails: mme4@georgetown.edu, i.m.karabash@gmail.com
}

\vspace{4ex}

\begin{abstract}
We derive Euler-Lagrange equations for the topology optimization of decay rate in 3-d lossy optical  cavities. This leads to a new  class of time-harmonic differential or integro-differential equations, which can be written as nonlinear Maxwell systems with switching functions of special types. Our approach is based on the notion of Pareto optimal frontier and on the multi-parameter perturbation theory for eigenfrequencies. Parallels with optimal control theory are discussed.
\end{abstract}

\vspace{1ex}


{
Resonators\\[0.5ex]
MSC-classes: 
35Q61, 
35B34, 
78M50, 
49R05, 
58E17, 
90C29, 
47B44, 
47J10, 
47A55 
\\[0.5ex]


\tableofcontents
}

\section{Introduction}

Optical cavities (or resonators) with high quality factor (Q-factor) are important components in contemporary optical engineering with a number of applications \cite{V03,SYZWKGX16,SKTN16,SDR18} ranging from quantum electrodynamics \cite{HR06,VLMS02} to bio-imaging \cite{S11}. Depending on a particular purpose, various designs are employed including Fabry-P\'{e}rot cavities composed of mirrors or micromirrors of special shape \cite{HR06,SDR18,KBP19}, toroidal resonators supporting whispering gallery modes \cite{V03,SYZWKGX16}, and photonic crystal cavities based on photonic bandgap effect \cite{AASN03,NKT08,SKTN16}.

Many applications and, in particular, quantum computing, require strong localization of light in a small volume. Since  radiation loss is, roughly speaking, in inverse proportion to the resonator size, the task of fabrication of microresonators with  high Q-factors is very nontrivial. A substantial attention was attracted by 2-d and 1-d photonic crystal designs because the parallel development of nanotechnology allowed researchers to fabricate sophisticatedly tailored  cavities  \cite{AASN03,NKT08,DMTSH14} and apply numerical methods of structural (or topology) optimization \cite{LSV03,BPSZ11,LJ13,MSG17,AFKRZ18}.

The numerical topology optimization may be based either on parametrization of certain class of designs by a handful of variables and running over them optimization algorithms    \cite{AASN03,NKT08,DMTSH14,MSG17}, or, in the case of the ``full topology optimization'', on the consideration of every ``pixel''
of a discretized version of the resonator as a degree of freedom \cite{LJ13}. 

While microcavities with very high Q-factor have been fabricated, the ratio between predicted by simulations and actually observed Q-factors reaches one or even two orders of magnitude (see \cite{DMTSH14} and references in \cite{MSG17}). This essential difference was attributed  to the neglection of absorption in simulations and to fabrication errors \cite{MSG17}.

While further progress may be achieved by the study of robustness of high-Q designs and by optimization in presence of uncertainties (see \cite{DMTSH14}), an essential obstacle is that 
the optimization of Q-factor in 3-d optical cavities is not  understood well enough from mathematical point of view. 

The aims of this note are theoretical study of the full topology optimization for eigenmode Q-factors in 3-d optical cavities and the derivation of an analogue of the Euler-Lagrange equation for this problem.  

While the employed cost functions appear to be standard and are connected with exponential loss rate (and Q-factor), they are used 
in the settings of Pareto optimization \cite{BV04,K13,AK17,KKV18}, which is more common in mathematical economics, than in Physics. These settings resolve the difficulties with ill-posedness of the Q-factor maximization (the issues with nonexistence of optimizers were discussed in detail in \cite{LJ13}). 
Our approach shifts the difficulty to the sensitivity analysis, which requires in the case of Pareto optimization a more advanced multi-parameter perturbation theory \cite{K14}.

\subsection{\label{ss:i} Eigenproblem and Pareto optimization settings} 

We employ the eigenproblem-based formulation, which involves the time-harmonic \linebreak 
Maxwell system
 \begin{align} 
   \ii \nabla \times \H (\x) & =  \omega \epr (\x) \ep_0 \E (\x), \label{e:E}
   \\
 - \ii \nabla \times \E  (\x) & =   \omega \mu_0 \H (\x) \label{e:H}
\end{align}
equipped with lossy boundary conditions responsible to the leakage of energy. The leakage leads to the appearance of eigenfrequencies $\om$ in the lower complex half-plane $\CC_- := \{ z \in \CC : \Im z <0\}$. The Silver-Müller radiation condition, absorbing layers, or impedance  boundary conditions  can be used for the modeling \cite{ACL17,LJ13,M03} of the dissipation effects. 

In the optimization process, the spatially dependent (relative  electric) permittivity $\epr (\x)$ is modified inside a fixed bounded domain $D_\c$ assuming that the constraints
\begin{gather} \label{e:c}
0< \ep_- \le \epr (\x) \le \ep_+ , \quad \x \in D_\c , 
\end{gather}
are satisfied.  Here $\ep_-$ and $\ep_+$ are minimal and maximal  permittivities of the materials available for the fabrication,
for example, $\ep_-=1$ (air) and $\ep_+ =11.9716 $ (silicon).
The permittivity functions $\epr (\cdot)$ satisfying  (\ref{e:c}) will be called feasible.

The problems of the maximization of the Q-factor and the minimization of the loss rate
have been intensively investigated experimentally \cite{AASN03,KTTKRN10}, theoretically \cite{VLMS02,K13,LJ13,OW13,KLV17,AFKRZ18,KKV18}, and numerically \cite{LSV03,HBKW08,KS08,NKT08,OW13,LJ13,MSG17,LFGBKBIWSHL18}. The quality-factor associated with a particular eigenmode $[\E \ \H ]$ of (\ref{e:E}), (\ref{e:H}) can be defined as $Q = \frac{|\Re \om|}{-2 \Im \om}$, where the real part $\al = \Re \om$ of $\om$ is the real-frequency of 
the corresponding eigenoscillations $\exp (-\ii \om t) [ \E (x) \ \H(x) ] $, and the negative of imaginary part $\Ga (\om) = - \Im \om$ is the exponential loss rate. In the case of radiation conditions,
$\Ga (\om)$ is the  half-bandwidth of a resonance $\om$.

The theoretical study of the maximization of the Q-factor requires a careful choice of an objective and of optimization constraints because many ``straightforward'' formulations of the problem are ill-posed \cite{LJ13}. The sequence of reformulations  suggested in  \cite{LJ13} to resolve this difficulty employs a number of ideas including the minimization of the volume of cavity for a given quality factor,
the single complex-frequency formulation, 
and the frequency-averaged  local density of states as objectives. In the somewhat different context of a quantum cavity, the issue with nonexistence of physically interesting minimizers was noticed already in the pioneering paper \cite{HS86}, where additional constraints $\al_1 \le \Re \om \le \al_2$ were suggested as one of possible recipes to avoid this difficulty. 

The goal of present paper is to return to the simple loss rate  objective $\Ga (\om)$, but to impose extreme versions of the real-frequency and cavity volume constraints with the aim to resolve the ill-posedness issue. That is, we assume that both a certain real-frequency $\Re \om = \al$ and the region $D_\c$, where the optimization is performed, are fixed. We take the loss rate $\Ga (\om) = - \Im \om$ as the cost function. So, for a given $\al>0$, $\Ga_{\min} (\al)$ is defined as the infimum of  all possible $\Ga$ such that $\om = \al - \ii \Ga$ is an achievable complex-frequency, i.e., an eigenfrequency generated by a certain feasible function $\epr$. 

Denoting the best possible value of eigenfrequency by $\om_{\min} (\al) = \al - \ii \Ga_{\min} (\al)$, we can formulate
the optimization problem in the following way:
\begin{multline*} 
\text{find feasible permittivity functions $\epr $} \\
\text{ that generate $\om_{\min} (\al)$ as their  eigenfrequency.} 
\end{multline*}
We call such $\epr (\cdot)$ optimal. In other words, we are interested in the maximal value of Q-factor for a prescribed value $\al$ of real-frequency $\Re \om$.

For the 1-d equation $y'' (s)= - \epr(s) \ep_0 \mu_0 \om^2 y (s)$, where $s=x_3$
and where $y= E_2$ describes a TEM-mode in a layered cavity, the similar formulation was theoretically and computationally studied in \cite{K13,KLV17,KKV18}. 
Generally, it is difficult to guarantee that the real parts $\Re \om $ of  achievable $\om $ cover the whole positive semi-line $\RR_+ = \{\al \in \RR: \al >0 \}$. Such real parts $\Re \om$ will be called achievable and their set is denoted by  $\dom \Ga_{\min} $. The 1d results of \cite{KLV17,KKV18} suggest that all high enough real parts $\Re \om$ are achievable.

The above theoretical formulation is natural from the point of view of the Pareto optimization \cite{BV04}. Indeed, the curve 
$\{ \al - \ii \Ga_{\min} (\al) : \al \in \dom \Ga_{\min}\} $ is the Pareto optimal frontier for achievable complex-eigenfrequencies $\om$ \cite{K13,K14}.

\subsection{\label{ss:ir}Main results and methods of the paper}
 
Our main goal is to derive the nonlinear Maxwell-type equation
\begin{align} 
& \frac{|\Phi_\E (\x)|}{\ep_0 \mu_0} \, \nabla \times \nabla \times \E (\x)  =  \om^2  \left( \frac{\ep_+ +\ep_-}{2} |\Phi_\E (\x)| +\frac{\ep_+-\ep_-}{2} \Phi_\E (\x) \right)  \E (\x), \label{e:ELh} 
\end{align} 
where  $\Phi_\E (\x) $ plays the same role as  \emph{switching functions}  do in optimal control theory. 
In the case of 3d full topology Pareto optimization the switching function can be written in the form 
\[
\Phi_\E (\x)
=  \Im (\E (\x) \cdot \E (\x)) , \quad \text{where 
$ \E \cdot \E = \E^\top  \E  = \sum_{j=1}^3 (E_j)^2$} 
\]
and $\E= [ E_1  \ E_2  \ E_3 ] = ( E_1 \ E_2  \ E_3 )^\top $ is perceived as a $\CC^3$-column vector. This equation is fulfilled 
inside the optimization region $D_\c$ for the $\E$-component of one of eigenmodes
$[\E  \ \H ] $ of (\ref{e:E}), (\ref{e:H}) corresponding to the optimal $\om = \om_{\min} (\al)$. 

We will say that (\ref{e:ELh}) equipped with a leaky boundary condition is an \emph{Euler-Lagrange eigenproblem} (EL-eigenproblem).  
In Section \ref{s:Imp},  the EL-eigenproblem is derived in the simplest possible settings, where $\om$ is assumed to be a  simple  complex-eigenfrequency of the Maxwell system with a boundary condition of impedance type. Then, in Section \ref{s:Abs}, we show that the form of the EL-eigenproblem does not substantially depend on a particular choice of the type of lossy boundary conditions (while optimal eigenmodes and permittivity functions $\epr$ obviously do depend). To this end a more realistic model of a leaky cavity with an absorbing layer in the outer domain is used. In Section \ref{s:Dis} we consider the connection between optimizers $\epr (\x)$ and the eigenfields  $\E (\x)$ of the EL-eigenproblems, and discuss some of related mathematical questions, in particular, the parallels with optimal control theory.

Since 3d photonic crystals with sophisticatedly tailored design are difficult in fabrication, the most essential efforts in the experimental and numerical studies were concentrated on structures with 2d and  1d geometries \cite{AASN03,HBKW08,NKT08,K13,KLV17,LFGBKBIWSHL18,KKV18}.

Photonic crystals with 2d (or 1d) structure described in $D_\c$ 
by $\epr (\x)$ depending only on 2  components (or a 1 component) of the $\x$-vector, but surrounded by a 3d outer medium, are considered in Section \ref{s:2d1d}. Using the same tools as in the 3d-geometry case, we show that corresponding EL-equations  have the form (\ref{e:ELh}) with different switching functions $\Phi_\E $ involving integrals over cross-sections of the optimization domain $D_\c$, which is assumed now to be  cylindrical. 
For example, in the case of optimization of a 2d photonic crystal we obtain 
the switching function 
\[
\Phi_\E (x) = \int_{c_-}^{c_+} \Im (\E (\x) \cdot \E (\x)) \dd x_3 ,
\]
where $\{ \x=(x_1,x_2,x_3) \in \RR^3 \ : \ c_-< x_3<c_+ \}$ is the slab containing the optimization domain $D_\c$.
So, in the case of 1d and 2d geometries, equation (\ref{e:ELh}) becomes of an integro-differential type.

The main tool for the derivation of the EL-eigenproblem is an analytical lemma about the sets covered by 2-parameter perturbations of zeros of analytic functions, which is a particular case of a results of \cite{K14} and is given in Appendix for a convenience of the reader. 

To apply the 2-parameter lemma, we need as a preliminary step a formula for the first correction term of an eigenfrequency of (\ref{e:E})-(\ref{e:H}) under 1-parameter perturbations of $\epr (\x)$. This perturbation formula is derived in Section \ref{ss:Per}.
While there is an abundance of studies on perturbation expansions for various types of Maxwell equations (see \cite{AFKRZ18,JJWM08,K13,KLV17,LLYBH90} and references therein), we were not able to find in the literature a corresponding perturbation result  for a lossy cavity. For resonances corresponding to outgoing conditions, a somewhat  less explicit formula can be found in \cite{LLYBH90}, where the degenerate case is also included. 
On the other side, the formula of Section \ref{ss:Per} can be seen as an extension to the case of 3d Maxwell equations of the  perturbation results  on the 3d Schrödinger equation  \cite{HS86} and on the 1d string case  \cite{K13,KLV17}. In 1d settings, the leading correction term of Puiseux  series for resonances of arbitrary multiplicity was studied analytically in \cite{K13,KLV17}.

\section{Euler-Lagrange eigenproblem in the case of \\
impedance boundary conditions}
\label{s:Imp}

Let us consider a bounded 3d domain $D$ (i.e., a bounded open connected subset of $\RR^3$) with a smooth boundary $\pa D$ and impose on $\pa D$ the impedance boundary condition of the simplest form 
\begin{equation}\label{e:IBC}
 \n \times \E - \Z \H_\tau= 0, \qquad \mbox{ for } \x \in \pa D ,
\end{equation}
where $\n$ is the unit outward normal vector to $\pa D$, $\H_\tau = (\n\times \H) \times \n$ is the tangential component of the field
$\H$, and the impedance coefficient $\Z (\x)$ is a scalar uniformly positive $L^\infty (\pa D)$-function (uniformly positive means that there exists a constant $C>0$ such that $\Z (\x) \ge C$ on the domain of its definition).

An optimization domain $D_\c$ is a 3d subdomain of $D$ that represents the region  where the scalar relative permittivity function $\epr (\cdot)$ can be changed in the process of optimization. The structure of the cavity  in the outer region $D_{\out} = D \setminus D_\c$ is assumed to be known, unchangeable during the optimization, however, not necessarily homogeneous. 
We assume that this outer structure is described by a uniformly positive permittivity function $\epout (\x)$, $\x \in D_{\out}$.

The family of all $L^\infty (D)$-functions $\vep_r $ such that  
\begin{align}
\text{$\ep_- \le \epr (\x) \le \ep_+$ for $\x \in D_\c$, \quad and $\vep_r (\x) = \vep_{out} (\x)$   for  $\x \in D_{\out}$,} \label{e:F}
\end{align}
will be called the family of feasible permittivities and denoted by $F$.

Assume that $\al>0$, that the optimal complex-eigenfrequency $\om_0 = \om_{\min} (\al) $ (see Section \ref{ss:i}) is achieved for a certain  feasible permittivity function 
$\vep_r (\x)= \vep_{\min} (\x) \in F$, and that for this $ \vep_{\min} $ the eigenfrequency $\om_0$ is 
 simple (i.e. of algebraic multiplicity 1). Let $\Psi_0 (\x) = [\E^0 (\x) \ \H^0 (\x)]$, $\x \in D$, be a  corresponding eigenfield, which is unique up to a multiplication on a complex constant.

We consider linear perturbations $\zeta p +\frac{1}{\vep_{\min} }$ in the direction $p(\x)$, $\x \in D$, of the inverse  $1/\vep_{\min} (\x)$ to the permittivity function. A direction $p $, which is also an $L^\infty (D)$-function,
 is called admissible 
if $\vep_r (\x,\zeta) = \left( (\vep_{\min} (\x))^{-1}  + \zeta p (\x) \right)^{-1}$ is feasible for small enough positive numbers $\zeta$.  
The perturbed eigenfrequency $\om (\zeta)$ corresponding to $\vep_r (\x,\zeta) $ is an analytic complex valued function in $\zeta$ for complex $\zeta$ close enough to $0$ and satisfies 
\begin{align} \label{e:omz}
\om (\zeta) & = \om_0 + \zeta \C_1 (p)  + O (|\zeta|^2) \text{ as } \zeta \to 0, 
\\
\text{where } \quad 
\C_1 (p) & = \om_0 \K_1 
 \int_{D_\c}  p   \vep_{\min}^2   \E^0 \cdot \E^0 \ \dd \x 
  \label{e:C1}
\end{align}
and $ \K_1  =  \left( \int_D \left( \vep_{\min}  \ep_0 \E^0  \cdot \E^0  - \mu_0 \H^0  \cdot \H^0  \right) \dd \x \right)^{-1} \neq 0  
$
(see Section \ref{ss:Per}, where this formula is analytically derived in more general settings).

Note that any admissible direction $p$ satisfies $p(\x) = 0$ in $D_\out$, and so, $p$ can be considered as $L^\infty (D_\c)$-function.
The first correction term $\C_1 $ depends also on $\om_0$, $\vep_{\min}$, and $\Psi_0$, but since they are fixed this dependence is not important for our needs. Then the mapping $p \mapsto \C_1 (p)$ can be considered as a complex-valued linear bounded functional on $L^\infty (D_\c)$. 

For the convenience in the use of Convex Analysis techniques, we identify the complex plane $\CC$ with the 2d real plane $\RR^2$. So $p \mapsto \C_1 (p)$  will be simultaneously considered as an $\RR^2$-valued and as a $\CC$-valued real-linear mapping defined on  $L^\infty (D_\c)$. 

 Since the family $F$ of feasible permittivities is closed and convex in $L^\infty (D)$, the set of admissible directions 
 \[
 A = \{ C (1/\vep -1/\vep_{\min}) \ : \ \vep \in F, \ C \ge 0 \}
 \]
  and its image 
 $\C_1 [A] = \{ \C_1 (p) \ :  \ p \in A \}$ under the mapping $p \mapsto \C_1 (p)$ are closed convex cones in $L^\infty (D_\c)$ and $\RR^2$, respectively.
 
Obviously, the origin $0$ belongs to the image $\C_1 [A]$ because one can take $p \equiv 0$ in $D_\c$. Moreover, 
\begin{gather*} \label{e:0inB}
\text{$0$ is on the boundary $\pa \C_1 [A]$ of $\C_1 [A]$.}
\end{gather*}
Otherwise, one can take two admissible directions $p_1 (\x)$ and $p_2 (\x)$ so that the open triangular domain 
\[
T = \{\C_1 (\zeta_1 p_1 + \zeta_2 p_2) \ : \ \zeta_1>0, \  \zeta_2>0, \  \zeta_1 +\zeta_2 <1 \},
\]
which is the interior of the convex hull of $0$, $\C_1 (p_1)$, and $\C_1 (p_2)$, contains a certain complex segment $\ii (0,\de) = \{\ii \tau \ : \ 0 < \tau < \de \} $ of the positive imaginary semi-axis $\ii \RR_+$. Then the Maclaurin series (\ref{e:omz}) and the homotopy arguments (see Appendix \ref{s:2Per}) easily lead to the conclusion that $\om_0 - \ii \de_1$ is an achievable eigenfrequency for a certain $\de_1>0 $, and so, to the contradiction with the assumption that $\om_0$ has a minimal possible loss rate $\Ga (\om)$ among all achievable $\om$ on the line 
$\{ \Re \om = \al \}$.

Since $\C_1 [A]$ is a closed convex cone having $0$ on its boundary, there exists $\th \in \RR$ such that the closed half-plane $e^{\ii \th} \overline{\CC}_+ = \{ e^{\ii \th} z : \Im z \ge 0\}$ contains $\C_1 [A]$.

Summarizing, we see from the form (\ref{e:C1}) of $\C_1 (p)$ that replacing, if necessary, $\Psi_0$  by its  multiple  $[\E \ \H ] = e^{\ii \th_1} \Psi_0$ with a certain phase $\th_1 \in \RR $ one can ensure that
\[
\text{  $\left\{ \textstyle \int_D  p \vep_{\min}^2  \E \cdot \E \ \dd \x  \ : \ p \in A \right\}$ is a subset of the upper closed half-plane  $ \overline{\CC}_+ $}.
\]
This and the constraints $\ep_- \le \epr (\x) \le \ep_+$ of the feasible family $F$ give for $\x \in D_\c$ that 
$\vep_{\min} (\x) = \ep_\pm$ whenever $\pm \Im (\E (\x) \cdot \E (\x)) >0$.

Combining the last observation with the Maxwell system, one gets for $\x \in D_\c$
\begin{equation} 
\ii \, |\Im (\E \cdot \E)| \ \nabla \times \H  =  \om_0 \ep_0  \left( \frac{\ep_+ +\ep_-}{2} |\Im (\E \cdot \E)| +\frac{\ep_+-\ep_-}{2} \Im (\E \cdot \E)\right) \E\label{e:ELE} 
\end{equation}
and, in turn, the EL-eigenproblem (\ref{e:ELh}) with $\Phi_\E = \Im (\E \cdot \E)$ and $\om = \om_0$.

\begin{rem} 
One can choose other forms of the switching function. For example, 
$\Phi_E = \Re (\E \cdot \E)$ leads to essentially equivalent EL-eigenproblem with a solution $\E $ replaced by $e^{-\ii \pi/4} \E $.
Note that the EL-equation is linear with respect to (w.r.t) multiplication of $\E $ on real constants $C$, but is nonlinear w.r.t. multiplication on complex constants $C \in \CC \setminus \RR$.
\end{rem}

\section{Model with absorbing conductivity layers}
\label{s:Abs}

A cavity in Optical Engineering is usually a part of a bigger complex structure. The leakage of energy occurs rather due to various dissipative effects in the medium near the cavity, than due to escape of waves to the far surrounding modeled by an idealized infinite vacuum. One can make the model of a lossy cavity more realistic and adjustable to various specific setting by the introduction of artificial absorbing region.  

On the other hand, absorption layers are intensively employed in  simulations because of their convenience from the point of view of numerical approximations \cite{M03}. 

A simplest way to model the absorption  is to introduce a spatially varying positive conductivity $\si (\x)$. This approach also allows one to take into account possible metallic and semiconductor structures in  the surrounding region $D_\out$.

For a time being, let us impose minimal assumptions on $\si (\x)$
supposing that it is a nonnegative $L^\infty (D)$-function. 
This model leads to the time-harmonic Maxwell system 
 \begin{align} 
   \ii \nabla \times \H (\x) - \ii \si (\x) \E (\x)& =  \omega \epr (\x) \ep_0 \E (\x)  \label{e:Esi}
   \\
 - \ii \nabla \times \E  (\x) & =   \omega \mu_0 \H (\x)  \label{e:Hsi}.
\end{align}
As before we assume that the impedance boundary condition (\ref{e:IBC}) is imposed on $\pa D$.

We will denote the corresponding pseudo-Hamiltonian \cite{E12} by 
$\M_{\epr^{-1}, 1,\si} $, where 
\[
\M_{\vep^{-1}, \mu^{-1},\si} \begin{pmatrix} \E \\ \H \end{pmatrix}  =  
\begin{pmatrix}
- \frac{\ii \si(\x)}{\vep (\x) \ep_0}  & \frac{\ii}{\vep (\x) \ep_0}  \nabla \times \\  -\frac{\ii}{\mu (\x) \mu_0} \nabla \times & 0 
\end{pmatrix} 
\begin{pmatrix} \E \\ \H \end{pmatrix} 
\]
is defined analogously  to classical conservative Hamiltonians
$\M_{\vep^{-1},\mu^{-1},0}$ of \cite{DL,JJWM08} 
 with the only difference that, instead of conservative boundary conditions on $\pa D$, the dissipative impedance boundary condition participates in the definition of the domain of the operator \cite{ELN,LL}.

\subsection{One-parameter perturbations of eigenfrequencies}
\label{ss:Per}

In this subsection we derive the first correction term for perturbations of 
a simple isolated eigenvalue $\om_0 $ of the linear operator 
$\MM_0 = \M_{\epr^{-1}, 1, \si}$. Let 
$ \Psi_0 = [ \E^0 \ \H^0 ]$ be an associated eigenfield. 

In the case of lossy cavity, $\MM_0$ is not Hermitian  and $\om_0$ is not necessarily real. 
Therefore an additional step is needed in the derivation of the first correction term. 
The resulting formula does not exactly coincide with the well-known first-order correction for a lossless cavity \cite{JJWM08}.  
This additional step requires  the use of the operator $\MM_0^*$ adjoint to $\MM_0$ w.r.t. 
the energy form 
\begin{multline}
\< \Psi,\Psi \>_{\epr} = \int_D \left(\ep_0 \< \epr (\x)  \E (\x), \E (\x) \>_{\CC^3} + \mu_0 \<  \H (\x), \H (\x) \>_{\CC^3} \right)\dd \x = \\ =
\int_D \left(\ep_0 \epr (\x) \E (\x) \cdot \overline{\E (\x)} + \mu_0 \H (\x) \cdot \overline{\H (\x)} \right)\dd \x
,
\end{multline}
where $\<  \cdot , \cdot \>_{\CC^3}$ is a standard sesquilinear inner product in the 3d complex vector space $\CC^3$, and 
$\overline{\E} = [ \overline{E_1} \ \overline{E_2} \ \overline{E_3} ] $, 
$ \overline{\H} = [ \overline{H_1} \ \overline{H_2} \ \overline{H_3} ]$ consist of complex conjugates to the components of $\E$ and $\H$.

We use the observation that the field $\Psi_\star$ defined by 
$ \Psi_\star = \begin{pmatrix}  \overline{\E_0 } \\ - \overline{\H_0 } \end{pmatrix} $ is an adjoint eigenstate of $\MM_0^*$ in the sense that
\[
\MM_0^* \Psi_\star = \overline{\om_0} \ \Psi_\star .
\]

To derive the first correction term to $\om_0$ under the perturbations 
\[
\epr (\x,\zeta) = \left( \epr (\x)^{-1}  + \zeta p (\x) \right)^{-1}
\]
 of the permittivity function, we introduce  the 
operator-valued function 
$
\MM (\zeta) = \MM_0 + \zeta \MM_1,$  where $\MM_1 = \M_{p,0,0}$
and $\zeta$ is a complex number in a vicinity of $0$.

Since $\MM (\zeta)$ is analytic and of Kato's type (A) \cite[Sect. VII.2]{K13}, the perturbed eigenvalue $\Om (\zeta)$ is analytic near $\zeta = 0$ and  
$
\om (\zeta) = \om_0 + \C_1 \zeta + o (\zeta) .
$
Moreover, there exists an eigenstate $\Psi (\x,\zeta)= \begin{pmatrix} \E (\x,\zeta)  \\ \H (\x,\zeta) \end{pmatrix}$ analytically dependent on $\zeta$ and satisfying 
\[
\MM(\zeta) \Psi (\zeta) =  \om (\zeta)  \Psi (\zeta) , \qquad  \Psi (\zeta) = \Psi_0 + \zeta \Psi_1 + o(\zeta) .
\]

From $\MM(\zeta) \Psi (\zeta) =  \om (\zeta)  \Psi (\zeta)$, we see that
\[
(\MM_0 + \zeta \MM_1) (\Psi_0 + \zeta \Psi_1 + o(\zeta)) = 
(\om_0 + \zeta \C_1 + o (\zeta) ) (\Psi_0 + \zeta \Psi_1 + o(\zeta)).
\]
and so 
\begin{align*}
 \C_1 \Psi_0 & =  \MM_1 \Psi_0 + (\MM_0 - \om_0) \Psi_1  ,\\
\C_1 \<  \Psi_0 , \Psi_\star \>_{\epr} & =  \< \MM_1 \Psi_0 , \Psi_\star \>_{\epr}  + 
\< \Psi_1,  (\MM_0^* - \overline{\om_0}) \Psi_\star \>_{\epr}  =
\< \MM_1 \Psi_0 , \Psi_\star \>_{\epr} .
\end{align*}
Due to the assumption that $\om_0$ is a simple eigenvalue,
one has
\[
0 \neq \< \Psi_0 , \Psi_\star \>_{\epr}  = \int_D \left(\ep_0 \epr \E^0 \cdot \E^0 - \mu_0 \H^0 \cdot \H^0 \right) \dd \x 
\]
and  
\begin{gather}
\C_1 = \frac{ \<  \M_{p,0,0} \Psi_0 , \Psi_\star \>_{\epr} }{ \<  \Psi_0 , \Psi_\star \>_{\epr}} =
\frac{\int_D  \<  \epr  \ii p \ep_0^{-1} \nabla \times \H^0 , \overline{\E^0 } \>_{\CC^3} \dd \x}
{\< \Psi_0 , \Psi_\star\>_{\epr} }  
\notag \\  =
\frac{\int_D  p  \< \epr (\om_0 \epr + \ii \si \ep_0^{-1})\E^0 , \overline{\E^0} \>_{\CC^3} \dd \x}{\< \Psi_0 , \Psi_\star \>_{\epr} }
=
\frac{\int_D p (\om_0 \epr^2 + \ii \si \epr \ep_0^{-1})  (\E^0 \cdot \E^0) \dd \x}{\< \Psi_0 , \Psi_\star \>_{\epr} }.
\label{e:C1si1}
\end{gather}
We see that the first correction term
$ \frac{\dd \om (\xi)}{\dd \xi} = \C_1 (p) $ is a linear functional of the perturbation direction $p$.


We apply (\ref{e:C1si1}) under additional assumption that $p(\x) \si(\x) = 0 $ everywhere in $D$ (cf. (\ref{e:si=0})), which leads to a simpler formula 
\begin{gather} \label{e:C1si}
\C_1 (p) = \om_0 \frac{\int_D p (\x) \epr^2 (\x) \sum_{j=1}^3 (E^0_j (\x))^2  \dd \x}{\int_D \left(\ep_0 \epr (\x) \sum_{j=1}^3 (E^0_j (\x))^2   - \mu_0 \sum_{j=1}^3 (H^0_j (\x))^2  \right) \dd \x}.
\end{gather}

\begin{rem}
For the lossless cavity, (\ref{e:C1}) can be reduced to the standard formula \cite{JJWM08}  since $\om$ is real and so the the real and imaginary parts of the eigenfields $\E^0$ and $\H^0$ are also  eigenfields. 
\end{rem}

\subsection{EL-equations with absorbing layers and local optimizers}
\label{ss:Abs}

The procedure of derivation of EL-equation in the presence of conductivity is similar to Section \ref{s:Imp}. However, we extend the EL-eigenproblem also to the case of local extrema.

Assume that the function $\epr $ belongs to the family of feasible permittivities $F$. That is, 
it is fixed in $D_\out$ and is a subject of optimization inside the optimization domain $D_\c$  satisfying the constraints (\ref{e:F}). 
Suppose that $\si (\cdot)$ is fixed during the optimization in whole $D$ and additionally
\begin{align} \label{e:si=0}
\text{$\si (\x) = 0$ for all $\x \in D_\c$,}
\end{align}
i.e., the medium is lossless inside the optimization domain $D_\c$.
So the loss of energy happens in a subset of $D_\out$ where $\si (\x) >0$ 
and on the boundary $\pa D$ due to the impedance boundary condition.

Let $\al>0$ and $\Re \om_0 = \alpha$. 
Assume that $\om_0$ is a simple eigenvalue of the Maxwell system (\ref{e:Esi}), (\ref{e:Hsi}), (\ref{e:IBC}) with a certain $\epr  \in F$.

To define local minimizers, let us introduce the family $F (\alpha)$ of 
pairs $\{ \om , \epr \}$ that satisfy the following conditions:
(i) $\Re \om = \al$, (ii) $\om$ is an eigenvalue of $M_{\epr^{-1},1,\si}$,
(iii) $\epr (\x)$ is a feasible permittivity function. 
That is, $F (\alpha)$ is a family of feasible pairs $\{ \om , \epr \}$ under the additional constraint $\Re \om =\alpha$, and can be considered as a subset of the normed space $\CC \times L^\infty (D)$ with any reasonable norm $\| \cdot \|$ (e.g., $\| \{ \om , \epr \} \| = \max \{ |\om|, \| \epr \|_{L^\infty} \}$).

Assume that $\{ \om_0 , \epr \}$ is a local minimizer in the sense that $\Ga (\om_0) \le \Ga (\om)$ for all pairs $\{ \om, \vep \} \in F (\al)$
in a ball of the space $\CC \times L^\infty (D)$ centered in $\{ \om_0 , \epr \}$ with a certain positive radius. Then applying the perturbation formula (\ref{e:C1si}) and the arguments of Section \ref{s:Imp}, one sees that there exists an eigenfield $[ \E \ \H ] $ of (\ref{e:Esi}), (\ref{e:Hsi}), (\ref{e:IBC}) that satisfies the EL-equation (\ref{e:ELE}) in the optimization domain $D_\c$.

The form of EL-equation in $D_\c$ does not depend on the way how the leakage of energy is modeled in the outer region. However, the corresponding eigenfield $[ \E \ \H ] $ in $D_\c$ does generally depend and so does the local optimal permittivity $\epr$.
Assuming that in the surrounding region $D_\out$ there exists an absorption subdomain $D_\abs$, where $\si(x) >0$ (and that $D_\abs$ has a positive volume),
one sees that the EL-eigenfield $[\E \ \H ]$ satisfies in $D_\out$ a system different 
from that of Section 2. The glue-type  condition  between the parts of eigenfield in $D_\c$ and $D_\out$ is encoded in the domain of definition of the pseudo-Hamiltonian $\M_{\epr^{-1},1,\si}$ (namely, in the requirements that the distributions $\nabla \times \E$ and 
$\nabla \times \H$ are $L^2 (D)$-vector-functions). 

To include this dependence on the absorption into the EL-system, let us introduce in the whole $D$ the constraint functions $\vep_\pm (\cdot)$ defined by 
\begin{equation*}
\text{$\vep_\pm (\x) = \ep_\pm $ for $\x \in D_\c$, and $\vep_\pm (\x) = \vep_\out (\x) $ 
for $\x \in D_\out$,} \label{e:eprpm}
\end{equation*}
and write the EL-equation in the whole $D$ as 
\begin{multline} 
\frac{|\Phi_\E|}{\ep_0 \mu_0} \ \nabla\times\nabla\times \E  - \frac{\ii \om \si (\x)}{\ep_0 } \E  = \om^2   \left( \frac{\vep_+ (\x) +\vep_- (\x)}{2} |\Phi_\E| +\frac{\vep_+ (\x) -\vep_- (\x)}{2} \Phi_\E \right) \E.\label{e:ELE3} 
\end{multline} 
with $\Phi_\E = \Im (\E \cdot \E)$.

\begin{rem} 
If we drop the assumption that $\si (\x) = 0$ in $D_\c$, but suppose instead that $\si (\x)$, $\x \in D$, is fixed in the process of optimization of $\epr$,
the resulting equation (\ref{e:ELE3}) is the same, but another switching function of the form $\Phi_\E (\x) = \Im \left( \left(\epr^2 + \frac{\ii \si \epr}{\ep_0 \om}\right) (\E \cdot \E) \right)$ has to be used.
However, if one would like to take into account small conductivity effects inside $D_\c$, a reasonable optimization problem should include also optimization of $\si (\x)$, $\x \in D_\c$, over a certain range of values. Such a model leads to switching effects not only for $\epr$, but also for $\si$. The resulting EL-equation is more technically involved and will be considered elsewhere. 
\end{rem}

Let $[ \E \ \H ] $ be a solution to the linear Maxwell system (\ref{e:Esi}), (\ref{e:Hsi}), (\ref{e:IBC}) with a certain feasible $\epr$ such that 
the eigenfield $\E$ satisfies additionally the EL-eigenproblem (\ref{e:ELE3}). We will say that such $\E (\x)$ is an EL-eigenfield for a nonlinear eigenvalue $\om$.

\begin{rem} 
 Generally, there is no reason to expect the uniqueness of optimal 
$\epr$ (see \cite{HS86,AK17,KKV18} for the discussion of uniqueness in the 1d case and in the case of a quantum cavity) and so there is no reason to expect the uniqueness of an EL-eigenfield.
\end{rem}

\subsection{Optimization of 2d and 1d photonic crystals embedded into a 3d outer medium}
\label{s:2d1d}

In this section, the outer medium described by $\vep_\out (\x)$ and $\si (\x)$ has as before a 3d geometry, as well as the outer  domain $D_\out$. 

The optimization domain of the cavity is taken to be a cylinder 
\[
D_\c = D_2 \times [x_{3,1},x_{3,2}] = \{ (x_1,x_2,x_3) \in \RR^3 \ : \ (x_1,x_2) \in D_2 \text{ and } c_- < x_3 < c_+  \},
\]
where $D_2$ is a bounded 2d domain and $c_\pm \in \RR$ are certain constatnts so that $c_-<c_+$.

The internal structure of the cavity in $D_\c$ is assumed to have  either 2d, or 1d geometry. To model this, suppose that the  permittivity function $\epr (\x)$ in the optimization domain $D_\c$ either depends only on $(x_1,x_2)$, or only on $x_3$. That is, $\epr (\x)$, $\x \in D_\c$, belongs
\begin{itemize}
\item either to the feasible 
family $F_2$ defined by the previous constraints (\ref{e:F}) and the additional assumption that for each $(x_1,x_2) \in D_2$ we have
$\epr (x_1,x_2,x_3) = \epr (x_1,x_2,x'_3) $ for all $x_3, x'_3 \in (c_-,c_+)$,

\item or, to the feasible family $F_1$ defined by  (\ref{e:F}) and the assumption that for $x_3 \in (c_-,c_+)$ we have
$\epr (x_1,x_2,x_3) = \epr (x'_1,x'_2,x_3) $ for all points 
$(x_1,x_2), (x'_1,x'_2) \in D_2$. 
\end{itemize}

Suppose that  $\si (\cdot)$ is fixed and $\si (\x) = 0$ for $\x \in D_\c$.

Similarly to the previous section, we assume that the pair $\{ \om , \epr \}$ is minimizer or local minimizer over $F_2$ or over $F_1$ for $\Ga (\om) = -\Im \om$  under the real part constraint $\Re \om_0 = \al>0$.

Then, the arguments of Section \ref{s:Imp} show that there exists 
an eigenfield $[ \E \ \H ]$ of the linear Maxwell system (\ref{e:Esi}), (\ref{e:Hsi}), (\ref{e:IBC}) that satisfies in $D_\c$
the EL-equation 
\begin{align*} 
& \frac{|\Phi_\E (\x)|}{\ep_0 \mu_0} \, \nabla \times \nabla \times \E (\x)  =  \om^2  \left( \frac{\ep_+ +\ep_-}{2} |\Phi_\E (\x)| +\frac{\ep_+-\ep_-}{2} \Phi_\E (\x) \right)  \E (\x), 
\end{align*} 
with  the switching function $\Phi_\E $ of the form 
\begin{itemize}
\item $\Phi_\E (\x) = \int_{c_-}^{c_+} \Im (\E (\x) \cdot \E (\x)) \dd x_3 $ in the case of optimization over the family $F_2$ of 2d photonic crystals;
\item 
$\Phi_\E (\x) = \int_{D_2} \Im (\E (\x) \cdot \E (\x)) \dd x_1 \dd x_2$ in the case of optimization over the family $F_1$ of 1d photonic crystals.
\end{itemize}
So the EL-equation has a non-local character of an integro-differential equation in $D_\c$. Similarly to Section \ref{ss:Abs}, it can be written as a nonlinear eigenvalue problem (\ref{e:ELE3}),  (\ref{e:IBC}) in the whole domain $D$.

\section{\label{s:Dis}Conclusions and discussion}

We have derived new Euler-Lagrange equations for the Pareto decay rate minimization in 3-d, 2-d, and 1-d photonic crystals embedded into a structured 3-d outer medium with dissipative effects. 

The efficiency of this approach to computations of the Pareto optimal frontier of optimal eigen-frequencies was up to now demonstrated only for the 1-d case of TEM-modes in idealized layered cavities \cite{KLV17,KKV18}. Application of the same scheme in 3-d settings require a deeper study of a series of new mathematical questions, which are briefly discussed below. 

Let us consider the dependence of locally optimal $\epr$ on the EL-eigenfield $\E $ in more detail. 
Section \ref{s:Imp} shows that, inside $D_\c$, we have $\epr (\x) = \ep_\pm$ whenever $\pm \Phi_\E (\x) >0$. However, (\ref{e:ELE}) does not impose any assumption on $\epr (\x) $ in the set $S$ of points $\x \in D_\c$ where $\Phi_\E (\x) =0 $.
Such phenomenon is common for some optimal control problems,
for example, for economic trading model and the moon landing problem \cite{SL12}. Adapting optimal control terminology, it is natural to call $S$ the  \emph{singular set} for the switching function $\Phi_\E$ (for connections between 1d resonance optimization and optimal control see \cite{KKV18}).

Since $S$ is the zero-level set of the switching function, it is natural to expect that, generically,
it is a 2d manifold and so has a zero volume (i.e., zero 3d Lebesgue measure). In such a case, the values of  $\epr$ on $S$ can be discarded since any change of $L^\infty$-coefficient 
$\epr $ on a set of zero volume does not influence the solution of the Maxwell system. Then adjusting Optimal Control terminology, it is natural to say that optimal $\epr (\x) $ has a \emph{bang-bang structure} in the sense that 
it takes in $D_\c$ only two extreme possible values $\ep_\pm$. From the Engineering point of view, this corresponds to the high contrast design, where 
air holes are sophisticatedly tailored in a silicon crystal.

However, it is difficult to rule out the existence of 3d settings where a switching function has a singular sets of positive volume. The main reason to expect that there exist pathological settings of such type is the fact that presently available unique continuation results require Lipschitz regularity of the coefficient $\epr$ \cite{V,HL}. In the optimization problem under consideration $\epr$ has, generally, a lower regularity of the $L^\infty$-space.

\begin{rem}
 Since the dynamic electromagnetic field can be considered to be real-valued, the complexification of the fields is due solely to the time-harmonic setting. One observes that $\Phi_\E = \Im (\E \cdot \E)= 2 \Re \E \cdot \Im \E$. So the equality $\Im (\E \cdot \E) = 0$
means that $\Re \E = [\Re E_1 \ \Re E_2 \ \Re E_3] $ and $ \Im \E = [\Im E_1 \ \Im E_2 \ \Im E_3]$ are orthogonal on $S$.
\end{rem}

\appendix

\section{Appendix: 2-parameter perturbations of eigenvalues}
\label{s:2Per}

Let us consider a function  $Q(z;\zeta)= Q(z;\zeta_1,\zeta_2)$ of 3 complex variables $z$, $\zeta_1$, $\zeta_2$ assuming that this function is analytic in a polydisc $\DD_R^3$, where $\DD_R = \{ \om \in \CC \ : \ |\om| <R\}$.
The variables $\zeta_1$ and $\zeta_2$ are united in the pair $\zeta = (\zeta_1,\zeta_2) \in \CC^2$ and play the role of perturbation parameters. We are interested in the behavior of the set of zeros 
of the function $Q(\cdot;\zeta)$. That is, we are interested in the set $\Sigma_Q (\zeta) $ of complex numbers $z \in \DD_R$ satisfying the equation $Q (z;\zeta) = 0$, and in the local evolution of this set under small perturbations of $\zeta$.

We take real positive solid triangles 
$ 
T_\de  :=  \{ \zeta \in \DD_R^2 \ : \  \zeta_1,
\zeta_2 > 0  \ \text{ and } \ 0 < \zeta_1 + \zeta_2 < \delta \} 
$ 
as families of feasible perturbation parameters. Then $\Sigma_Q [T_\de]  := \bigcup_{\zeta \in T_\de} \Sigma_Q (\zeta) $ is the set of $z$-zeros achievable over $T_\de$.

Assume that 
\begin{itemize}
\item 
$z_0 = 0$ is a simple zero of a function $Q (\cdot, 0,0)$, 
\item  
$\eta_j := -\frac{\pa Q (0,0,0)}{\pa \zeta_j}\left(\frac{\pa Q (0,0,0)}{\pa z} \right)^{-1} \neq 0 , \quad j=1,2,
$
\item and $\th_j = \arg \eta_j$ are ordered such that 
$\th_2 = \th_1 + \th_0 \ (\modn 2\pi) $ with a certain $\th_0 \in (0,\pi)$.
\end{itemize}

Then there exists a small enough $r>0$ and an analytic function $\om (\zeta)$ defined for $\zeta \in \DD_r^2$ such that
$
Q (\om (\zeta); \zeta) = 0 
$ everywhere in $\DD_r^2$ and 
\begin{gather} \label{e:Ts}
\om (\zeta) = \om_1 (\zeta)  + o (\zeta) \quad \text{ with  $ \om_1 (\zeta) = \eta_1 \zeta_1 + \eta_2 \zeta_2$ as $\zeta \to 0$.}
\end{gather}
It can be seen from (\ref{e:Ts}) that for triangles $T_\delta$ with small enough $\delta$
the images
$\om_1 [T_\de] = \{ \om_1 (\zeta) \ : \zeta \in T_\delta \}$ and 
$\om [T_\de] = \{ \om (\zeta) \ : \zeta \in T_\delta \}$ lie inside a certain complex half-plane, for example, inside 
$\exp \left( \ii (\th_1 - \de_0) \right) \CC_+$ with $\de_0$ small enough (we use here the assumptions  that $0<\th_0 < \pi$ and $\th_1 < \th_2 <\th_1 + \pi$).

Now, considering the evolution as $\de \to 0$ of images $\Si_Q [\pa T_\de]$ of the triangle's boundaries $\pa T_\de$,  it is easy to show the main statement of this section:
for any small enough $\de_n>0$, $n=1,2$, there exists 
 $\de_3>0$ such that
\[
\{ z \in \DD_{\de_3} \ : \ \th_1 + \de_2 < \arg z < \th_2 - \de_2 \} \subset \Sigma_Q [ T_{\de_1} ]  .
\]

The above arguments are the simplest particular case of a more general result \cite{K14}  on multi-parameter perturbations of zeros of arbitrary finite multiplicity. 
The proof in \cite{K14} is more technically involved because Puiseux series have to be used for multiple zeros instead of the Taylor series (\ref{e:Ts}).

\vspace{1ex}
\noindent
\emph{Acknowledgments.} 
ME and IK have been supported by the VolkswagenStiftung project “Modeling, Analysis, and Approximation Theory toward applications in tomography and inverse problems”.
IK has been supported by the Alexander von Humboldt Foundation.

\end{document}